\documentclass[11pt,a4paper]{article}
\usepackage[margin=1in,includefoot]{geometry}
\usepackage{amssymb}
\usepackage{amsmath}
\usepackage{bm}
\usepackage{array}
\usepackage{float}
\usepackage{stmaryrd}
\usepackage{multirow}
\usepackage{amsthm}
\usepackage{graphicx,multicol,caption,lipsum}
\usepackage{subcaption}
\usepackage{xcolor}
\usepackage[shortlabels]{enumitem}
\usepackage{mathtools}
\usepackage[colorlinks,linkcolor=blue,citecolor=blue]{hyperref}
\usepackage[noabbrev]{cleveref}
\usepackage{tikz}
\usepackage{pgfplots} 
\usepgfplotslibrary{polar}
\pgfplotsset{compat=1.5}
\usepackage{soul}
\usepackage{authblk}

\title{Least squares spectral element formulation of eigenvalue problems with/without interface : the one dimensional example}
\author[$\dag$]{Himanshu Garg}
\author[$\ddag$]{Fleurianne Bertrand}
\author[$\dag$]{Subhashree Mohapatra\thanks{Email id: subhashree@iiitd.ac.in(Corresponding author)}}
\affil[$\dag$]{Department of Mathematics, IIIT Delhi, India}
\affil[$\ddag$]{Department of Mathematics, TU Chemnitz, Germany}
\date{March 2025}
\begin{document}
\date{~}
\maketitle
\begin{abstract}
Here, we present a least-squares based spectral element formulation for one-dimensional eigenvalue problems with interface conditions. First we develop the method for without interface case, then we extend it to interface case. Convergence analysis for eigenvalues and eigenfunctions have been discussed. Numerical experiments with different jump conditions have been displayed.
    \end{abstract}

    \newtheorem{thm}{Theorem}[section]
    \newtheorem{prop}{Proposition}[section]
    \newtheorem{lem}{Lemma}[section]
    \newtheorem{coro}{Corollary}[section]
    \newtheorem{rem}{Remark}[section]
    \newtheorem{guess8}{Example}[section]

		\newcolumntype{C}{>{\centering\arraybackslash} m{2cm} }
		\newcolumntype{D}{>{\centering\arraybackslash} m{1cm} }
		\newcolumntype{E}{>{\centering\arraybackslash} m{2.4cm} }
         \newcolumntype{F}{>{\centering\arraybackslash} m{2.8cm} }
         \newcolumntype{G}{>{\centering\arraybackslash} m{2.8cm} }
         \newcolumntype{H}{>{\centering\arraybackslash} m{2.8cm} }
\section{Introduction}  
One-dimensional eigenvalue problems play an important role in various physics and engineering related problems as they can be used to investigate various complex systems after reducing them to a single variable. For example, these kind of problems can be used to understand the behavior of waves and vibrations, states of a quantum system, modes of a classical mechanical oscillators, natural frequencies of vibrations in structures, stability of various control systems.

Pierce and Varga \cite{pi1972} have presented high order approximations of linear eigenvalue problems in several dimensions. Schrodinger equation using a finite element formulation has been solved in \cite{ja}. A finite difference based formulation for one-dimensional Schrodinger equation has been discussed in \cite{tr}. One dimensional eigenvalue problems on unbounded domains are presented in \cite{ho}.

We discuss a least-squares based spectral element formulation for one-dimensional eigenvalue problems with interface conditions. Interface conditions are weakly imposed in the least-squares functional, which will lead to discrete variational formulation, that is a generalized eigenvalue problem in finite dimensional setup that is preconditioned with the preconditioner obtained using the equivalent quadratic forms which is then solved using the \emph{eigs} command of MATLAB.

In Section 2, we discuss variational eigenvalue problems and their well posedness. Discretization and stability estimates have been presented in Section 3. Discrete eigenvalue problems have been discussed in Section 4. Convergence analyis of eigenvalues and eigenfunctions have been discussed in Section 5. Results on numerical experiments have been displayed in Section 6. Finally we conclude with Section 7.
\section{Mathematical background}
In this section, we introduce eigenvalue problems for second-order ordinary differential equations with and without interface conditions in their strong forms. Their equivalent weak forms are derived under suitable regularity assumptions. We discuss the existence of eigenpairs of the weak problems using theory of compact-self adjoint operators.
\subsection*{Without interface}
Here, we consider the eigenvalue problem for the second order differential equation of the form
\begin{equation}
    \mathcal{L}y(x)=-\left(p(x)y'(x)\right)'+q(x)y(x)=\lambda r(x)y(x),\quad x\in\Omega=(a,b)\label{egvp_prob}
\end{equation}
subject to certain boundary conditions, labeled as \eqref{dir_cd}-\eqref{periodic_cd} below. These boundary conditions include the Dirichlet, Neumann, and periodic cases.
\begin{align}
    y(a)=0,\quad\quad &y(b)=0;\label{dir_cd}\tag{I}\\
    y'(a)=0,\quad\quad &y'(b)=0;\label{neu_cd}\tag{II}\\
    y(a)=y(b),\quad\quad &y'(a)=y'(b).\label{periodic_cd}\tag{III}
\end{align}
The problem of solving the equation \eqref{egvp_prob} subject to the boundary condition \eqref{dir_cd} will be referred to as problem (I) and similar notations will be used for other problems.

Now, we make some assumptions on the coefficients that $p\in C^{2}(\overline{\Omega})$ and $q,r\in C^{1}(\overline{\Omega})$. Also, we assume that there exists $p_{0}>0$ such that $$p(x)\geq p_{0}>0\quad\text{and}\quad q(x),r(x)>0\quad\forall x\in\overline{\Omega}.$$
In case of periodic boundary conditions, we make a further assumption on $p$ such that $p(a)=p(b).$

Now, consider the Hilbert spaces $V$, where $V=H^{1}_{0}(\Omega)=\{v\in H^{1}(\Omega):\hspace{2pt}v(a)=v(b)=0\}$ for Dirichlet case, $V=H^{1}(\Omega)$ for Neumann case, and $V=H^{1}_{p}(\Omega)=\{v\in H^{1}(\Omega):\hspace{2pt}v(a)=v(b)\}$ for periodic boundary conditions case, with norm $\|\cdot\|_{1,\Omega}$ and $V_{1}=L^{2}(\Omega)$ with norm $\|\cdot\|_{\Omega}$. Then, the variational eigenvalue problem for problems \eqref{dir_cd}-\eqref{periodic_cd} is: Find an eigenpair $(\lambda,u)$ such that
\begin{equation}
    a(u,v)=\lambda b(u,v),\quad\quad\forall v\in V,\label{var_egvp}
\end{equation}
where
$a(\cdot,\cdot)$ is a symmetric bilinear form on $V$ which is also elliptic and continuous and $b(\cdot,\cdot)$ is a continuous, symmetric and positive definite bilinear form on $V_{1}$ and are defined by
\begin{align*}
    a(u,v)&=\int_{a}^{b}\left[p(x)u'(x)v'(x)+q(x)u(x)v(x)\right]dx,\quad\quad\forall u,v\in V,\\
    b(u,v)&=\int_{a}^{b}r(x)u(x)v(x)dx,\quad\quad\forall u,v\in V_{1}.
\end{align*}
Consider the associated source problem: Given $f\in V_{1}$, find $u\in V$ such that 
\begin{equation}
    a(u,v)=b(f,v),\quad\quad \forall v\in V.\label{src_prblm}
\end{equation}
Since, $a(.,.)$ is a continuous, elliptic bilinear form in $V$ and $b(f,\cdot)$ is a bounded linear functional on $V$ so by Lax-Milgram theorem $\exists!$ $T:V_{1}\rightarrow V$, the corresponding solution operator , such that $Tf=u$ and 
\begin{equation}
    a(u,v)=a(Tf,v)=b(f,v),\quad\forall v\in V.\label{var_src_prblm}
\end{equation}
Also,
\begin{equation*}
    \|u\|_{1,\Omega}=\|Tf\|_{1,\Omega}\leq C\|f\|_{\Omega}.
\end{equation*}
Since $V$ is compactly embedded in $V_{1}$ by the Rellich theorem, $T:V\rightarrow V$ is a compact operator. Moreover, $T:V\rightarrow V$ is also a self-adjoint operator since the bilinear forms $a(\cdot,\cdot)$ and $b(\cdot,\cdot)$ are symmetric. Then using the spectral theorem for self-adjoint compact operators, we get that there exists countably infinite many eigenvalues  $\{\mu_{i}\}_{i=1}^{\infty}\subset\mathbb{R}$ of $T$ which converges to $0$. Now, it can be seen that $(\mu,u)\in\mathbb{R}\times V$ is an eigenpair for $T$ if and only if $\left(\frac{1}{\mu},u\right)$ is an eigenpair of \eqref{var_egvp}. Therefore, there are countably infinite many eigenvalues of \eqref{var_egvp} which will diverge to $\infty$ and can be written as
\begin{equation}
    0<\lambda_{1}<\lambda_{2}<\lambda_{3}<\dots
\end{equation}Henceforth, all the analysis will be done on the operator $T$ to derive the results for the variational eigenvalue problem.

Now to develop our numerical scheme, we further assume the higher regularity for our solution operator that is 
\begin{equation}
    \|Tf\|_{2,\Omega}\leq C\|f\|_{\Omega}.\label{reg_est}
\end{equation}
\subsection*{With interface}
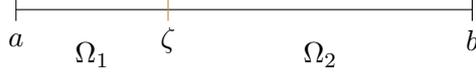
\begin{figure}[ht]
	\centering
	\begin{tikzpicture}
		\draw (0,0) -- (6,0);
		\draw (0,-0.15) --(0,0.15);
		\draw (6,-0.15) --(6,0.15);
		\draw[color=brown] (2,-0.15) -- (2,0.15);
		\node at (0,-0.4) {$a$};
		\node at (6,-0.4) {$b$};
		\node at (2,-0.4) {$\zeta$};
		\node at (1,-0.6) {$\Omega_{1}$};
		\node at (4,-0.6) {$\Omega_{2}$};
	\end{tikzpicture}
\caption{One dimensional interface problem }
\label{fig:intf}
\end{figure}

We consider the following eigenvalue problem for second-order  ordinary differential equation \eqref{egvp_intf_prob} having an interface at $x=\zeta\in\Omega=(a,b)$, which partition the domain $\Omega$ into $\Omega_{1}=(a,\zeta)$ and $\Omega_{2}=(\zeta,b)$ as shown in \autoref{fig:intf},  subject to Dirichlet B.C. and also the interface conditions \eqref{intf_jump} and \eqref{intf_der_jump} mentioned below:
\begin{align}
    \mathcal{L}u=-\frac{d}{dx}\left(p\frac{du}{dx}\right)+q u&=\lambda ru,\quad\text{in }\Omega_{1}\cup\Omega_{2}\label{egvp_intf_prob}
\end{align}
with 
\addtocounter{equation}{-1}
\begin{subequations}
    \begin{align}
        \left\llbracket u\right\rrbracket&=0\label{intf_jump}\\
    \left\llbracket p\frac{du}{dx}\right\rrbracket&=0\label{intf_der_jump}\\
    u(a)=0 \quad u(b)&=0
    \end{align}
\end{subequations}
where $\llbracket\cdot\rrbracket$ denotes the jump of the function across the interface $\zeta$, defined as
$$\left\llbracket u\right\rrbracket:=u(\zeta-)-u(\zeta+),\quad\text{and}\quad \left\llbracket p\frac{du}{dx}\right\rrbracket=p(\zeta-)\frac{du}{dx}(\zeta-)-p(\zeta+)\frac{du}{dx}(\zeta+)$$
and the coefficient $p$ is piecewise $C^{2}$ function. For simplicity we write 
\begin{align*}
    p(x)=\begin{cases}
        p_{1}(x),& x\in(a,\zeta),\\
        p_{2}(x),& x\in(\zeta,b).
    \end{cases}
\end{align*}
where $p_{1}\in C^{2}(\overline{\Omega}_{1})$ and $p_{2}\in C^{2}(\overline{\Omega}_{2})$ and $q,r\in C^{1}(\overline{\Omega}).$ Also, we assume that there exists $p_{0}>0$ such that
$$p(x)\geq p_{0}>0,\forall x\in\Omega_{1}\cup\Omega_{2}\quad\quad q(x),r(x)>0,\quad\forall x\in\overline{\Omega}.$$

Now, consider the Hilbert space $V=H^{1}_{0}(\Omega):=\{v\in H^{1}(\Omega):\hspace{2pt} v(a)=v(b)=0\}$ with norm $\|\cdot\|_{1,\Omega}$ and $V_{1}=L^{2}(\Omega)$ with norm $\|\cdot\|_{\Omega}$. Then, the variational eigenvalue problem for the interface problem is: Find an eigenpair $(\lambda,u)$ such that
\begin{equation}
    a(u,v)=\lambda b(u,v),\quad\forall v\in V.\label{intf_weak_egvp}
\end{equation}
where $a(\cdot,\cdot)$ and $b(\cdot,\cdot)$ are same as of without interface case.

Consider the associated source problem: Given $f\in V_{1}$, find $u\in V$ such that
\begin{equation}
    a(u,v)=b(f,v),\quad\quad \forall v\in V.\label{intf_weak_src}
\end{equation}
Applying the Lax-Milgram theorem $\exists$ $T_{intf}:V_{1}\rightarrow V$, the corresponding solution operator, such that $T_{intf}f=u$ and 
\begin{equation}
    a(u,v)=a(T_{intf}f,v)=b(f,v),\quad\quad v\in V.
\end{equation}
Also, with same reasoning there are countably infinite many eigenvalues of \eqref{intf_weak_egvp} which will diverges to $\infty$ and cane be written as
$$0<\tilde{\lambda}_{1}<\tilde{\lambda}_{2}<\tilde{\lambda}_{3}<\dots$$
where $(1/\tilde{\lambda}_{i})$ is an eigenvalue of $T_{intf}$ solution operator.

Also, we assume the further that the weak solution has a higher regularity and the estimate
$$\|Tf\|_{\tilde{H}^{2}(\Omega)}\leq \|f\|_{\Omega}$$
holds true.
\section{Discretization and stability estimates}
In this section, we discuss the domain discretization and stability estimates related to spectral element functions. First, we discuss for eigenvalue problems without interface conditions and then we consider eigenvalue problems with interface case.
\subsection{Discretization}
\subsubsection*{Without interface}
\begin{figure}[ht]
	\centering
	\begin{tikzpicture}
		\draw (0,0) -- (6,0);
		\draw (0,-0.15) --(0,0.15);
		\draw (6,-0.15) --(6,0.15);
		\draw (1,-0.1) -- (1,0.1);
		\draw (2,-0.1) -- (2,0.1);
		\draw (3,-0.1) -- (3,0.1);
		\draw (4,-0.1) -- (4,0.1);
		\draw (5,-0.1) -- (5,0.1);
		\node at (0,-0.4) {$x_{0}=a$};
		\node at (6,-0.4) {$x_{L}=b$};
            \node at (1,-0.4) {$x_{1}$};
		\node at (2,-0.4) {$x_{2}$};
            \node at (3,-0.4) {$x_{3}$};
            \node at (5,-0.4) {$x_{L-1}$};
            \node at (4,-0.4) {\dots};
		\node at (3,-0.8) {$\Omega$};
	\end{tikzpicture}
\caption{Discretization of domain $\Omega$}
\label{domain_disc}
\end{figure}
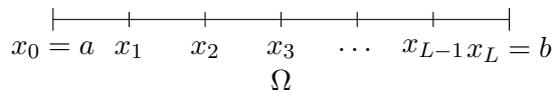

For discretization, let $a=x_{0}<x_{1}<x_{2}<\dots<x_{L}=b$ be the partition of the domain $\Omega$ into $L$ number of spectral elements $\{\Omega_{l}\}_{l=1}^{L}$ as shown in \autoref{domain_disc} such that for each $l=1,\dots,L$, we have $\Omega_{l}=(x_{l-1},x_{l})$ and $h_{l}=x_{l}-x_{l-1}$ be the length of that element. Also, define $h=\max_{1\leq l\leq L}\{h_{l}\}$ be the global mesh size of the partition.

Now for each $l=1,\dots,L$, define an invertible linear map $M_{l}:\overline{\hat{\Omega}}\rightarrow\overline{\Omega}_{l}$, where $\hat{\Omega}=(-1,1)$ is the master element, defined by
$$M_{l}(\xi)=x_{l-1}\left(\frac{1-\xi}{2}\right)+x_{l}\left(\frac{1+\xi}{2}\right).$$
Following the construction of elements, let us define the spectral element function $\{\hat{u}_{1},\hat{u}_{2},\dots,\hat{u}_{L}\}$ on $\hat{\Omega}$ to be the polynomials of degree at most $W$, that is for each $l\in\{1,\dots,L\}$ 
$$\hat{u}_{l}(\xi)=\sum_{i=0}^{W}c_{i}^{l}\xi^{i},$$
for some $c_{i}^{l}\in\mathbb{R}$.
Using these we define the space $\Pi^{W}$ to be the collection of all spectral element functions, that is
$$\Pi^{W}:=\left\{ \{\hat{u}_{l}\}_{l=1}^{L}:\hspace{2pt}\hat{u}_{l}\in P_{W}(\hat{\Omega}),\hspace{2pt}\forall l=1,\dots,L\right\},$$
where $P_{W}(\hat{\Omega})$ is the space of polynomials of degree upto $W$ with real coefficients restricted to $\hat{\Omega}$. Also, the basis for $P_{W}(\hat{\Omega})$ is the Lagrange interpolating polynomials of degree $W$ defined over the Gauss-Lobatto Legendre (GLL) nodes. Let us define the norm on $\Pi^{W}$ such that
$$\|\{\hat{u}_{l}\}\|^{2}_{W}:=\sum_{l=1}^{L}\|\hat{u}_{l}\|^{2}_{2,\hat{\Omega}}.$$
Let us now define the set of non-conforming element functions over the physical elements to be the collection of $\{u_{l}\}_{l=1}^{L}$ where
$$u_{l}:\Omega_{l}\rightarrow\mathbb{R}\quad\text{such that} \quad u_{l}(x):=\hat{u}_{l}(M^{-1}_{l}(x)),\quad\forall x\in\Omega_{l}.$$
We now state the fundamental regularity result for function $u\in H^{2}(\Omega)$ based on \eqref{reg_est}, which is stated as follows:\\
Under the assumptions for $\mathcal{L}$ mentioned earlier, we have for $u\in H^{2}(\Omega)$
\begin{subequations}
\begin{align}
    \|u\|_{2,\Omega}&\leq C\{\|\mathcal{L}u\|_{\Omega}+|u(a)|+|u(b)|\}\quad\quad\text{(Dirichlet case)}\label{dirc_reg_est}\\
    \|u\|_{2,\Omega}&\leq C\{\|\mathcal{L}u\|_{\Omega}+|u'(a)|+|u'(b)|\}\quad\quad\text{(Neumann case)}\label{neu_reg_est}\\
    \|u\|_{2,\Omega}&\leq C\{\|\mathcal{L}u\|_{\Omega}+|u(a)-u(b)|+|u'(a)-u'(b)|\}\quad\quad\text{(Periodic B.C. case)}\label{periodic_reg_est}
\end{align}
\end{subequations}
\subsubsection*{With interface}
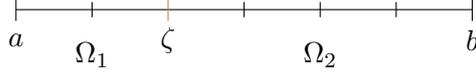
\begin{figure}
	\centering
	\begin{tikzpicture}
		\draw (0,0) -- (6,0);
		\draw (0,-0.15) --(0,0.15);
		\draw (6,-0.15) --(6,0.15);
		\draw (1,-0.1) -- (1,0.1);
		\draw[color=brown] (2,-0.15) -- (2,0.15);
		\draw (3,-0.1) -- (3,0.1);
		\draw (4,-0.1) -- (4,0.1);
		\draw (5,-0.1) -- (5,0.1);
		\node at (0,-0.4) {$a$};
		\node at (6,-0.4) {$b$};
		\node at (2,-0.4) {$\zeta$};
		\node at (1,-0.6) {$\Omega_{1}$};
		\node at (4,-0.6) {$\Omega_{2}$};
	\end{tikzpicture}
	\caption{Domain with interface discretization}
    \label{dom_disc_intf}
\end{figure}

To discretize the domain with interface, we partition each domain $\Omega_{i}$ into $N_{i}$ number of elements respectively. Let $a=x_{0}^{1}<x^{1}_{1}<\dots<x_{N_{1}}^{1}=\zeta$ and $\zeta=x_{0}^{2}<x^{1}_{2}<\dots<x_{N_{2}}^{2}=b$ be the partition of the domains $\Omega_{i}$ into $N_{i}$ number of spectral elements $\{\Omega_{l}^{1}\}_{l=1}^{N_{1}}$ and $\{\Omega_{l}^{2}\}_{l=1}^{N_{2}}$ respectively as shown in \autoref{dom_disc_intf} such that for each $l=1,\dots,N_{1}$ and $m=1,\dots,N_{2}$ we have $\Omega_{l}^{1}=(x^{1}_{l-1},x^{1}_{l})$ and $\Omega^{2}_{m}=(x^{2}_{m-1},x^{2}_{m})$.
Also, define the length of elements by $h_{l}^{1}=x_{l}^{1}-x^{1}_{l-1}$ and $h_{m}^{2}=x_{m}^{2}-x_{m-1}^{2}$ and $h=\max\{h_{l}^{1},h_{m}^{2}\}$ be the global mesh size for the partition.

In a similar manner as before, there exist invertible maps $M_{l}^{1}:\hat{\Omega}\rightarrow\Omega_{l}^{1}$ and $M_{m}^{2}:\hat{\Omega}\rightarrow\Omega_{m}^{2}$ for $l=1,\dots,N_{1}$ and $m=1,\dots, N_{2}$.

Now, we define the space $\Pi_{I}^{W}$ to be the collection of all spectral element functions, that is
$$\Pi^{W}_{I}:=\left\{\left\{\{\hat{u}_{l}^{1}\}_{l=1}^{N_{1}},\{\hat{u}_{m}^{2}\}_{m=1}^{N_{2}}\right\}:\hspace{2pt}, \hat{u}_{l}^{1}\in P_{W}(\Omega_{l}^{1}),\hat{u}_{m}^{2}\in P_{W}(\Omega_{m}^{2}),\hspace{2pt}l=1,\dots,N_{1},\hspace{2pt}m=1,\dots, N_{2}\right\}$$
Let us now, define the set of non-conforming element functions over the physical elements to be the collection of $\left\{\{{u}_{l}^{1}\}_{l=1}^{N_{1}},\{{u}_{m}^{2}\}_{m=1}^{N_{2}}\right\}$ as defined earlier.
We now state the fundamental regularity result for function $u\in H^{2}(\Omega)$ based on \eqref{reg_est}, which is stated as follows:\\
Under the assumptions for $\mathcal{L}$ mentioned earlier, we have for $u\in \tilde{H}^{2}(\Omega)$
\begin{align}
    \|u\|_{\tilde{H}^{2}(\Omega)}&\leq C\left\{\|\mathcal{L}u\|_{\Omega}+\left|\left\llbracket p\frac{du}{dx}\right\rrbracket\right|+|u(a)|+|u(b)|\right\}
\end{align}
\subsection{Quadratic forms and its stability estimates}

We will now define the quadratic form which will play a crucial role to define the numerical scheme for the variational eigenvalue problem.
\subsubsection*{Without interface}
Firstly, let us define the differential operator $\mathcal{L}_{l}$, for $l=1,\dots,L$, such that for $\phi\in H^{2}(\hat{\Omega})$, $(\mathcal{L}_{l}\phi):\hat{\Omega}\rightarrow\mathbb{R}$ is defined by
$$\mathcal{L}_{l}\phi(\xi)=(\mathcal{L}\phi(M_{l}(\xi)))\sqrt{J_{l}(\xi)},\quad\forall \xi\in\hat{\Omega},$$
where $J_{l}$ denotes the Jacobian of the map $M_{l}$. Then using this operator we define $\mathcal{L}^{a}_{l}$, a differential operator obtained by replacing the coefficients of $\mathcal{L}_{l}$ by their polynomial approximations of degree $W-1$, obtained by interpolating the coefficients over the GLL nodes on $\overline{\hat{\Omega}}$. Also, it can be shown that 
$$\sum_{l=1}^{L}\int_{\Omega_{l}}|\mathcal{L}u_{l}(x)|^{2}dx\leq C\sum_{l=1}^{L}\int_{-1}^{1}|\mathcal{L}_{l}^{a}\hat{u}_{l}(\xi)|^{2}d\xi+\varepsilon_{W}\|\{\hat{u}_{l}\}_{l=1}^{L}\|^{2}_{W},$$
where $C$ is a constant and $\varepsilon_{W}\rightarrow0$ as $W\rightarrow\infty$. In fact, $\varepsilon_{W}$ is spectrally small in $W$.

Since we follow a non-conforming approach without any constraints on space, therefore jumps in the functions and its derivatives along the inter-element nodes are required to be defined and added in the numerical scheme.

Let $x_{i}$ be the common node for the two elements $\Omega_{i}$ and $\Omega_{i+1}$. Clearly, $x_{i}$ is the image of $\xi=1$ under the map $M_{i}$ and is the image of $\xi=-1$ under the map $M_{i+1}$ and by chain rule, we have
\begin{align*}
    \frac{du_{i}}{dx}(x_{i})&=\frac{d\hat{u}_{i}}{d\xi}(1)\left.\frac{d}{dx}(M_{i}^{-1}(x))\right|_{x=x_{i}}=\frac{2}{h_{i}}\frac{d\hat{u}_{i}}{d\xi},\\
    \frac{du_{i+1}}{dx}(x_{i})&=\frac{d\hat{u}_{i+1}}{d\xi}(-1)\left.\frac{d}{dx}(M_{i+1}^{-1}(x))\right|_{x=x_{i}}=\frac{2}{h_{i+1}}\frac{d\hat{u}_{i+1}}{d\xi},
\end{align*}
Then we define the jump in spectral functions and its derivatives along the inter-element nodes $x_{i}$ as 
\begin{align*}
    \left[\hat{u}\right]_{i}&:=\hat{u}_{i+1}(-1)-\hat{u}_{i}(1),\\
    \left[\frac{d\hat{u}}{d\xi}\right]_{i}&:=\frac{du_{i+1}}{dx}(x_{i})-\frac{du_{i}}{d\xi}(x_{i})=\frac{2}{h_{i+1}}\frac{d\hat{u}_{i+1}}{d\xi}-\frac{2}{h_{i}}\frac{d\hat{u}_{i}}{d\xi}.
\end{align*}
Now, we can define the quadratic form for problems \eqref{dir_cd}-\eqref{periodic_cd} as $\mathcal{V}^{W}:\Pi^{W}:\mathbb{R}$ such that
\begin{align}
    \label{quad_form}\mathcal{V}^{W}\left(\hat{u}\right):=\sum_{l=1}^{L}\|\mathcal{L}^{a}_{l}\hat{u}_{l}\|^{2}_{\hat{\Omega}}+\sum_{i=1}^{L-1}\left(\left|\left[\hat{u}\right]_{i}\right|^{2}+\left|\left[\frac{d\hat{u}}{d\xi}\right]_{i}\right|^{2}\right)+\tilde{\mathcal{V}}^{W}\left(\hat{u}\right),\quad \forall \hat{u}=\{\hat{u}_{l}\}_{l=1}^{L}\in\Pi^{W}
\end{align}
where $\tilde{\mathcal{V}}^{W}:\Pi^{W}:\mathbb{R}$ corresponds to the boundary conditions and is defined as follows:
\addtocounter{equation}{-1}
\begin{subequations}
\begin{align}
    \tilde{\mathcal{V}}^{W}\left(\{\hat{u}_{l}\}_{l=1}^{L}\right):=\begin{cases}
        |\hat{u}_{1}(-1)|^{2}+|\hat{u}_{L}(1)|^{2},& \eqref{dir_cd}\\
        \left|\frac{2}{h_{1}}\frac{d\hat{u}_{1}}{d\xi}(-1)\right|^{2}+\left|\frac{2}{h_{L}}\frac{d\hat{u}_{L}}{d\xi}(1)\right|^{2},& \eqref{neu_cd}\\
        |\hat{u}_{1}(-1)|^{2}-\hat{u}_{L}(1)|^{2}+\left|\frac{2}{h_{1}}\frac{d\hat{u}_{1}}{d\xi}(-1)-\frac{2}{h_{L}}\frac{d\hat{u}_{L}}{d\xi}(1)\right|^{2},& \eqref{periodic_cd}
    \end{cases}
\end{align}
\end{subequations}
Let us now state a necessary result for proving the stability estimate and also the stability estimate, proved in \cite{shivangi2023}:
\begin{lem}
    Let $\{\hat{u}_{i}\}_{i=1}^{L}\in\Pi^{W}$. Then, there exists $\{\hat{v}_{i}\}_{i=1}^{L}\in\Pi^{3}$ such that $\hat{v}_{1}(-1)=\hat{v}_{L}(1)=0$ and the function $w:\Omega\rightarrow\mathbb{R}$ defined by $w|_{\Omega_{i}}=(\hat{u}_{i}+\hat{v}_{i})\circ M^{-1}_{i}$ is in $H^{2}(\Omega)$. Also, there exists a constant $C$ such that
    $$\sum_{l=1}^{L}\|\hat{v}_{l}\|^{2}_{2,\hat{\Omega}}\leq C\sum_{i=1}^{L-1}\left(\left|\left[\hat{u}\right]_{i}\right|^{2}+\left|\left[\frac{d\hat{u}}{d\xi}\right]_{i}\right|^{2}\right)$$
\end{lem}

\begin{thm}\label{stb_est}
    For $W$ large enough, there exists a constant $C,\tilde{C}>0$ such that
    $$\frac{1}{C}\sum_{l=1}^{n}\|\hat{u}_{l}\|^{2}_{2,\hat{\Omega}}\leq \mathcal{V}^{W}\left(\{\hat{u}_{l}\}_{l=1}^{L}\right)\leq \tilde{C} \sum_{l=1}^{n}\|\hat{u}_{l}\|^{2}_{2,\hat{\Omega}}.$$
\end{thm}
\begin{rem}
    Theorem \autoref{stb_est} implies that the quadratic form and the norm $\|\cdot\|_{W}$ on the space of spectral elements $\Pi^{W}$ are equivalent. Therefore, the quadratic form corresponding to the norm $\|\cdot\|_{W}$ can be used as preconditioner to improve the condition number of the numerical scheme.
\end{rem}
\subsubsection*{With interface}
Let us define the differential operators $\mathcal{L}_{1,l}$ and $\mathcal{L}_{2,m}$, for $l=1,\dots,N_{1}$ and $m=1,\dots,N_{2}$, such that for $\phi\in H^{2}(\hat{\Omega})$, we have $$\mathcal{L}_{1,l}\phi=(\mathcal{L}\phi(M_{l}^{1}(\xi)))\sqrt{J_{l}^{1}(\xi)},\quad\quad \mathcal{L}_{2,m}\phi=(\mathcal{L}\phi(M_{l}^{2}(\xi)))\sqrt{J_{l}^{2}(\xi)},\quad\forall\xi\in\hat{\Omega}.$$
where $J_{l}^{1}$ and $J_{m}^{2}$ are the Jacobian of the respective maps $M_{l}^{1}$ and $M_{m}^{2}$. Then in a similar manner we define $\mathcal{L}^{a}_{1,l}$ and $\mathcal{L}^{a}_{2,m}$ operators, whose all coefficients are polynomial approximations.
Now, we can define the quadratic form as $\mathcal{V}^{W}_{I}:\Pi^{W}_{I}\rightarrow\mathbb{R}$ such that
\begin{multline}
    \mathcal{V}_{I}^{W}\left(\hat{u}_{1},\hat{u}_{2}\right):=\sum_{i=1}^{2}\left[\sum_{j=1}^{N_{i}}\|\mathcal{L}^{a}_{i,j}\hat{u}^{i}_{j}\|^{2}_{\hat{\Omega}}+\sum_{j=1}^{N_{i}-1}\left(\left|\left[\hat{u}_{i}\right]_{j}\right|^{2}+\left|\left[\frac{d\hat{u}_{i}}{d\xi}\right]_{j}\right|^{2}\right)\right]+|\hat{u}^{1}_{1}(-1)|^{2}+|\hat{u}^{2}_{N_{2}}(1)|^{2}\\
    |\hat{u}_{N_{1}}^{1}(1)-\hat{u}^{2}_{1}(-1)|^{2}+\left|\frac{2}{h^{1}_{N_{1}}}p_{1}(M_{N-{1}}^{1}(1))\frac{d\hat{u}^{1}_{N_{1}}}{d\xi}(1)-\frac{2}{h^{2}_{1}}p_{2}(M^{2}_{1}(-1))\frac{d\hat{u}^{2}_{1}}{d\xi}(-1)\right|^{2}
\end{multline}
We now state a necessary result for stability estimate and the estimate from \cite{shivangi2023}.
\begin{lem}
    Let $\{\hat{u}_{1},\hat{u}_{2}\}\in \Pi^{W}_{I}$ where $\hat{u}_{1}=\{\hat{u}^{1}_{l}\}_{l=1}^{N_{1}}$ and $\hat{u}_{2}=\{\hat{u}^{2}_{m}\}_{m=1}^{N_{2}}$. Then there exists $\{\hat{v}_{1},\hat{v}_{2}\}\in\Pi^{W}_{I}$, where $\hat{v}_{1}=\{\hat{v}^{1}_{l}\}_{l=1}^{N_{1}}$ and $\hat{v}_{2}=\{\hat{v}^{2}_{m}\}_{m=1}^{N_{2}}$, such that $\hat{v}^{1}_{1}(-1)=0,\hat{v}^{1}_{N_{1}}(1)=0$ and $\hat{v}^{2}_{1}(-1)=0,\hat{v}^{2}_{N_{2}}(1)=0$ which will give $w^{1}_{i}=u_{i}^{1}+v^{1}_{i}\in H^{2}(\Omega_{1})$ and $w^{2}_{i}=u^{2}_{i}+v^{2}_{i}\in H^{2}(\Omega_{2})$. Also, there exists a constant $C$ such that
    $$\sum_{i=1}^{N_{1}}\|\hat{v}_{i}^{1}\|^{2}_{2,\hat{\Omega}}+\sum_{j=1}^{N_{2}}\|\hat{v}^{2}_{j}\|^{2}_{2,\hat{\Omega}}\leq C\sum_{i=1}^{2}\left[\sum_{j=1}^{N_{i}-1}\left(\left|\left[\hat{u}_{i}\right]_{j}\right|^{2}+\left|\left[\frac{d\hat{u}_{i}}{d\xi}\right]_{j}\right|^{2}\right)\right]$$
\end{lem}
\begin{thm}
    For large enough $W$, there exists $C,\tilde{C}>0$ such that 
    $$\frac{1}{C}\left(\sum_{i=1}^{N_{1}}\|\hat{u}_{i}^{1}\|^{2}_{2,\hat{\Omega}}+\sum_{j=1}^{N_{2}}\|\hat{u}^{2}_{j}\|^{2}_{2,\hat{\Omega}}\right)\leq C\mathcal{V}^{W}_{I}\left(\hat{u}_{1},\hat{u}_{2}\right)\leq \tilde{C}\left(\sum_{i=1}^{N_{1}}\|\hat{u}_{i}^{1}\|^{2}_{2,\hat{\Omega}}+\sum_{j=1}^{N_{2}}\|\hat{u}^{2}_{j}\|^{2}_{2,\hat{\Omega}}\right).$$\label{stab_est_intf}
\end{thm}

\section{Discrete variational eigenvalue formulation}
\subsection*{Without interface}
Firstly, we define the minimization of the least squares functional $\mathcal{R}^{W}:\Pi^{W}\rightarrow\mathbb{R}$, as a numerical scheme for the corresponding source problem \eqref{src_prblm}, defined by
\begin{equation}
    \mathcal{R}^{W}(\hat{u}):=\sum_{l=1}^{L}\|\mathcal{L}^{a}_{l}\hat{u}_{l}-\hat{r}_{l}\hat{F}_{l}\|^{2}_{\hat{\Omega}}+\sum_{i=1}^{L-1}\left(\left|\left[\hat{u}\right]_{i}\right|^{2}+\left|\left[\frac{d\hat{u}}{d\xi}\right]_{i}\right|^{2}\right)+\tilde{\mathcal{R}}^{W}\left(\hat{u}\right),\quad \forall \hat{u}=\{\hat{u}_{l}\}_{l=1}^{L}\in\Pi^{W}
\end{equation}
where 
\addtocounter{equation}{-1}
\begin{subequations}
\begin{align}
    \tilde{\mathcal{R}}^{W}\left(\{\hat{u}_{l}\}_{l=1}^{L}\right):=\begin{cases}
        |\hat{u}_{1}(-1)|^{2}+|\hat{u}_{L}(1)|^{2},& \eqref{dir_cd}\\
        \left|\frac{2}{h_{1}}\frac{d\hat{u}_{1}}{d\xi}(-1)\right|^{2}+\left|\frac{2}{h_{L}}\frac{d\hat{u}_{L}}{d\xi}(1)\right|^{2},& \eqref{neu_cd}\\
        |\hat{u}_{1}(-1)|^{2}-\hat{u}_{L}(1)|^{2}+\left|\frac{2}{h_{1}}\frac{d\hat{u}_{1}}{d\xi}(-1)-\frac{2}{h_{L}}\frac{d\hat{u}_{L}}{d\xi}(1)\right|^{2},& \eqref{periodic_cd}\\
    \end{cases}
\end{align}
\end{subequations}
and $\hat{r}_{l}$ is the polynomial approximation of $r(M_{l}(\xi))(J_{l}(\xi))^{1/4}$ and $\hat{F}_{l}$ is the polynomial approximation of $f(M_{l}(\xi))(J_{l}(\xi))^{1/4}$. Since, minimizing the above function is equivalent to solving its normal equation, which leads to the variational formulation for the corresponding source problem, that is to find $\hat{u}=\{\hat{u}_{l}\}_{l=1}^{L}\in \Pi^{W}$ such that
\begin{equation}
    a_{W}(\hat{u},\hat{v})=b_{W}(\hat{F},\hat{v}),\quad\quad\forall \hat{v}=\{\hat{v}_{l}\}_{l=1}^{L}\in\Pi^{W}\label{disc_src_problm}
\end{equation}
where $\hat{F}=\{\hat{F}_{l}\}_{l=1}^{L}$ and the bilinear forms are defined by
\addtocounter{equation}{-1}
\begin{subequations}
\begin{align}
    a_{W}(\hat{u},\hat{v})&:=\sum_{l=1}^{L}(\mathcal{L}^{a}_{l}\hat{u}_{l},\mathcal{L}^{a}_{l}\hat{u}_{l})_{\hat{\Omega}}+\sum_{i=1}^{L-1}\left(\left[\hat{u}\right]_{i}\left[\hat{v}\right]_{i}+\left[\frac{d\hat{u}}{d\xi}\right]_{i}\left[\frac{d\hat{v}}{d\xi}\right]_{i}\right)+\tilde{a}_{W}(\hat{u},\hat{v}),\\
    \tilde{a}_{W}(\hat{u},\hat{v})&:=\begin{cases}
        \hat{u}_{1}(-1)\hat{v}_{1}(-1)+\hat{u}_{L}(1)\hat{v}_{L}(1),& \eqref{dir_cd}\\
        \left(\frac{2}{h_{1}}\right)^{2}\frac{d\hat{u}_{1}}{d\xi}(-1)\frac{d\hat{v}_{1}}{d\xi}(-1)+\left(\frac{2}{h_{L}}\right)^{2}\frac{d\hat{u}_{L}}{d\xi}(1)\frac{d\hat{v}_{L}}{d\xi}(1), & \eqref{neu_cd}\\
        \left(\hat{u}_{1}(-1)-\hat{u}_{L}(1)\right)\left(\hat{v}_{1}(-1)-\hat{v}_{L}(1)\right)&\\
        \quad\quad+\left(\frac{2}{h_{1}}\frac{d\hat{u}_{1}}{d\xi}(-1)-\frac{2}{h_{L}}\frac{d\hat{u}_{L}}{d\xi}(1)\right)\left(\frac{2}{h_{1}}\frac{d\hat{v}_{1}}{d\xi}(-1)-\frac{2}{h_{L}}\frac{d\hat{v}_{L}}{d\xi}(1)\right),&\eqref{periodic_cd}
    \end{cases}
    \intertext{and}
    b_{W}(\hat{F},\hat{v})&:=\sum_{l=1}^{L}(\hat{r}_{l}\hat{F}_{l},\mathcal{L}^{a}_{l}\hat{v}_{l})_{\hat{\Omega}}.
\end{align}
\end{subequations}
Using this formulation, we define the discrete variational eigenvalue problem in the non-conforming setting as, to find $\hat{u}=\{\hat{u}_{l}\}_{l=1}^{l}\in\Pi^{W}$  and $\lambda^{W}$ such that
\begin{equation}
    a_{W}(\hat{u},\hat{v})=\lambda^{W}b_{W}(\hat{u},\hat{v}),\quad\quad\forall \hat{v}=\{\hat{v}_{l}\}_{l=1}^{L}\in\Pi^{W}.
\end{equation}
\subsection*{With interface}
Now, we define the minimization of following least squares functional $\mathcal{R}_{I}^{W}:\Pi^{W}_{I}\rightarrow\mathbb{R}$, as a numerical scheme for the corresponding source problem, defined by
\begin{multline}
    \mathcal{R}^{W}_{I}(\hat{u}_{1},\hat{u}_{2}):= \sum_{i=1}^{2}\left[\sum_{j=1}^{N_{i}}\|\mathcal{L}^{a}_{i,j}\hat{u}^{i}_{j}-\hat{r}^{i}_{j}\hat{F}^{i}_{j}\|^{2}_{\hat{\Omega}}+\sum_{j=1}^{N_{i}-1}\left(\left|\left[\hat{u}_{i}\right]_{j}\right|^{2}+\left|\left[\frac{d\hat{u}_{i}}{d\xi}\right]_{j}\right|^{2}\right)\right]+|\hat{u}^{1}_{1}(-1)|^{2}+|\hat{u}^{2}_{N_{2}}(1)|^{2}\\
    |\hat{u}_{N_{1}}^{1}(1)-\hat{u}^{2}_{1}(-1)|^{2}+\left|\frac{2}{h^{1}_{N_{1}}}p_{1}(M_{N-{1}}^{1}(1))\frac{d\hat{u}^{1}_{N_{1}}}{d\xi}(1)-\frac{2}{h^{2}_{1}}p_{2}(M^{2}_{1}(-1))\frac{d\hat{u}^{2}_{1}}{d\xi}(-1)\right|^{2}
\end{multline}
Since, minimization of the above function is equivalent to solving its normal equation, which leads to the variational formulation for the corresponding source problem, that is to find $\hat{u}=(\hat{u}_{1},\hat{u}_{2})=(\{\hat{u}^{1}_{i}\}_{i=1}^{N_{1}},=\{\hat{u}^{2}_{j}\}_{j=1}^{N_{2}})\in\Pi^{W}_{I}$ such that
\begin{equation}
    a_{W,I}(\hat{u},\hat{v})=b_{W,I}(\hat{F},\hat{v}),\quad\quad\forall \hat{v}\in \Pi^{W}_{I}.
\end{equation}
where the bilinear forms are defined by
\begin{multline}
    a_{W,I}(\hat{u},\hat{v}):=\sum_{i=1}^{2}\left[\sum_{j=1}^{N_{i}}(\mathcal{L}^{a}_{i,j}\hat{u}^{i}_{j},\mathcal{L}^{a}_{i,j}\hat{v}^{i}_{j})_{\hat{\Omega}}+\sum_{j=1}^{N_{i}-1}\left(\left[\hat{u}_{i}\right]_{j}\left[\hat{v}_{i}\right]_{j}+\left[\frac{d\hat{u}_{i}}{d\xi}\right]_{j}\left[\frac{d\hat{v}_{i}}{d\xi}\right]_{j}\right)\right]\\
    +\hat{u}_{1}^{1}(-1)\hat{v}^{1}_{1}(-1)+\hat{u}^{2}_{N_{2}}(1)\hat{v}^{2}_{N_{2}}(1)+\left(\hat{u}_{N_{1}}^{1}(1)-\hat{u}^{2}_{1}(-1)\right)\left(\hat{v}_{N_{1}}^{1}(1)-\hat{v}^{2}_{1}(-1)\right)\\
    +\left[\left(\frac{2}{h^{1}_{N_{1}}}p_{1}(M_{N-{1}}^{1}(1))\frac{d\hat{u}^{1}_{N_{1}}}{d\xi}(1)-\frac{2}{h^{2}_{1}}p_{2}(M^{2}_{1}(-1))\frac{d\hat{u}^{2}_{1}}{d\xi}(-1)\right)\right.\\
    \times\left.\left(\frac{2}{h^{1}_{N_{1}}}p_{1}(M_{N-{1}}^{1}(1))\frac{d\hat{v}^{1}_{N_{1}}}{d\xi}(1)-\frac{2}{h^{2}_{1}}p_{2}(M^{2}_{1}(-1))\frac{d\hat{v}^{2}_{1}}{d\xi}(-1)\right)\right],
\end{multline}
and
\begin{equation}
    b_{W,I}(\hat{F},\hat{v}):=\sum_{i=1}^{2}\sum_{j=1}^{N_{i}}(\hat{r}_{j}^{i}\hat{F}^{i}_{j},\hat{v}^{i}_{j})_{\hat{\Omega}}.
\end{equation}
Using this formulation, we define the discrete variational eigenvalue in the non-conforming setting as, to find $\hat{u}=(\hat{u}_{1},\hat{u}_{2})\in\Pi^{W}_{I}$ and $\lambda^{W}_{I}$ such that
\begin{equation}
    a_{W,I}(\hat{u},\hat{v})=\lambda b_{W,I}(\hat{u},\hat{v}),\quad\quad\forall\hat{v}\in\Pi^{W}_{I}.
\end{equation}
\section{Convergence analysis}
In this section, we discuss convergence analysis of eigenvalues and eigenfunctions in case of without interface and with interface case.
\subsection*{Without interface}
Theorem \ref{stb_est} and the Lax-Milgram theorem will ensure that for each given $f\in V_{1}=L^{2}(\Omega)$ there exists a unique $\hat{u}\in\Pi^{W}$, solution of the corresponding discrete source problem $\eqref{disc_src_problm}$. Now, we state the error between this discrete solution and the solution of the \eqref{src_prblm} which is proved in \cite{shivangi2023}
\begin{thm}
    Let $f\in V$, where $V$ is defined above. Assume $u=Tf$ is the solution of the corresponding source problem \eqref{src_prblm} and $\hat{U}_{i}(\xi)=u(M_{i}(\xi))$ for $\xi\in\hat{\Omega}$. Also, let $\hat{u}\in\Pi^{W}$ be the minimizer of $\mathcal{R}^{W}$ over $\Pi^{W}$. Then for $W\geq3$, the following estimate holds.
    \begin{equation}
        \sum_{i=1}^{L}\|\hat{u}_{i}-\hat{U}_{i}\|^{2}_{2,\hat{\Omega}}\leq CW^{-2}\sum_{i=1}^{L}\|\hat{U}_{i}\|_{3,\hat{\Omega}}^{2}.
    \end{equation}
\end{thm}
After obtaining the non-conforming solution we make a linear correction on each element $\Omega_{i}$ to make the non-conforming spectral functions $\{u_{i}\}_{i}$, defined by the transformation $M_{i}$ on $\hat{u}_{i}$, so that the corrected spectral element functions $\{u^{c}_{i}\}_{i}$ are conforming and belongs to $H^{1}(\Omega)$. Then we have the following estimate in $H^{1}$ norm
$$\|u-u^{c}\|_{1,\Omega}\leq CW^{-1}\|f\|_{1,\Omega}$$
Therefore, we define our discrete solution operator $T^{W}:V\rightarrow V$ to be $T^{W}(f)=u^{c}$ where $u^{c}$ is the conforming spectral element on $\Omega$. Hence, we have
$$\|(T-T^{W})(f)\|_{1,\Omega}\leq CW^{-1}\|f\|_{1,\Omega}$$
which implies that $T^{W}$ converges to $T$  uniformly in norm.
Therefore, using the results in the spectral approximation theory of linear operators in \cite{kato1995} we get that $\lambda^{W}\rightarrow\lambda$, where $\lambda$ is an eigenvalue of $T$.
\begin{rem}
    Assume that the coefficients of differential equations are analytic in $\overline{\Omega}$. Let $\lambda_{j}^{W}$ be an eigenvalue of $T^{W}$ converges to $\lambda$, and let $u_{j}^{W}$ be the corresponding eigenvectors. Then there exists an eigenvector $u$ of $T$ corresponding to $\lambda$ such that
    \begin{align*}
        \|u-u_{j}^{W}\|_{1,\Omega}\leq C_{1}e^{-b_{1}W}\\
        |\lambda-\lambda_{j}^{W}|\leq C_{1}e^{-b_{1}W}.
    \end{align*}
\end{rem}
\subsection*{With interface}
Theorem \eqref{stab_est_intf} and the Lax-Milgram theorem will ensure that for each $f\in V_{1}$, there exists a unique $\hat{u}\in\Pi^{W}_{I}$, solution of the corresponding discrete source problem. Now, we state the error between this discrete solution and the solution of \eqref{intf_weak_src} which is proved in \cite{shivangi2023}.
\begin{thm}
    Let $f\in V$, where $V$ is defined above. Assume $u=T_{intf}f$ is the solution of the corresponding source problem \eqref{intf_weak_src} and $\hat{U}^{i}_{j}(\xi)=u(M^{i}_{j}(\xi)$ for $\xi\in\hat{\Omega}$. Also, let $\hat{u}\in\Pi^{W}_{I}$ be the minimizer of $\mathcal{R}^{W}_{I}$ over $\Pi^{W}_{I}$. Then for $W\geq 3$, the following estimate holds:
    $$\sum_{i=1}^{2}\sum_{j=1}^{N_{1}}\|\hat{u}^{i}_{j}-\hat{U}^{i}_{j}\|^{2}_{2,\hat{\Omega}}\leq C W^{-2}\sum_{i=1}^{2}\sum_{j=1}^{N_{i}}\|\hat{U}_{j}^{i}\|^{2}_{3,\hat{\Omega}}.$$
\end{thm}
After obtaining the non-conforming solution we make a linear correction on each element $\Omega^{i}_{j}$ to make the non-conforming spectral functions $\{u^{i}_{j}\}_{i,j}$, defined by the transformation $M^{i}_{j}$ on $\hat{u}^{i}_{j}$, so that the corrected spectral element functions $\{\tilde{u}^{i}_{j}\}_{i}$ are conforming and belongs to $H^{1}(\Omega)$. Then we have the following estimate in $H^{1}$ norm
$$\|u-\tilde{u}\|_{1,\Omega}\leq CW^{-1}\|f\|_{1,\Omega}$$
Therefore, we define our discrete solution operator $T^{W}_{I}:V\rightarrow V$ to be $T^{W}_
{I}(f)=u^{c}$ where $\tilde{u}$ is the conforming spectral element on $\Omega$. Hence, we have
$$\|(T_{intf}-T^{W}_{I})(f)\|_{1,\Omega}\leq CW^{-1}\|f\|_{1,\Omega}$$
which implies that $T^{W}_{I}$ converges to $T_{intf}$  uniformly in norm.
Therefore, using the results in the spectral approximation theory of linear operators in \cite{kato1995} we get that $\lambda^{W}\rightarrow\lambda$, where $\lambda$ is an eigenvalue of $T_{intf}$.
\begin{rem}
    Assume that the coefficients of differential equations are analytic in $\overline{\Omega}_{1}$ and $\overline{\Omega}_{2}$. Let $\lambda_{j}^{W}$ be an eigenvalue of $T^{W}_{I}$ converges to $\lambda$, and let $u_{j}^{W}$ be the corresponding eigenvectors. Then there exists an eigenvector $u$ of $T_{intf}$ corresponding to $\lambda$ such that
    \begin{align*}
        \|u-u_{j}^{W}\|_{1,\Omega}\leq C_{1}e^{-b_{1}W}\\
        |\lambda-\lambda_{j}^{W}|\leq C_{1}e^{-b_{1}W}.
    \end{align*}
\end{rem}
\section{Numerical results}
In this section, we present our numerical experiments on one-dimensional eigenvalue problems with interface conditions.

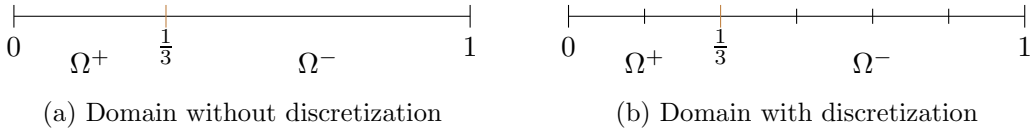
\begin{figure}[H]
\centering
\begin{subfigure}{0.45\textwidth}
	\centering
	\begin{tikzpicture}
		\draw (0,0) -- (6,0);
		\draw (0,-0.15) --(0,0.15);
		\draw (6,-0.15) --(6,0.15);
		\draw[color=brown] (2,-0.15) -- (2,0.15);
		\node at (0,-0.4) {$0$};
		\node at (6,-0.4) {$1$};
		\node at (2,-0.4) {$\frac{1}{3}$};
		\node at (1,-0.6) {$\Omega^{+}$};
		\node at (4,-0.6) {$\Omega^{-}$};
	\end{tikzpicture}
	\caption{Domain without discretization}
\end{subfigure}
\begin{subfigure}{0.45\textwidth}
	\centering
	\begin{tikzpicture}
		\draw (0,0) -- (6,0);
		\draw (0,-0.15) --(0,0.15);
		\draw (6,-0.15) --(6,0.15);
		\draw (1,-0.1) -- (1,0.1);
		\draw[color=brown] (2,-0.15) -- (2,0.15);
		\draw (3,-0.1) -- (3,0.1);
		\draw (4,-0.1) -- (4,0.1);
		\draw (5,-0.1) -- (5,0.1);
		\node at (0,-0.4) {$0$};
		\node at (6,-0.4) {$1$};
		\node at (2,-0.4) {$\frac{1}{3}$};
		\node at (1,-0.6) {$\Omega^{+}$};
		\node at (4,-0.6) {$\Omega^{-}$};
	\end{tikzpicture}
	\caption{Domain with discretization}
\end{subfigure}
\caption{One dimensional interface problem }
\label{int1d}
\end{figure}

We solve the following uni-dimensional boundary problem
\begin{align}\label{1d_int}
			-\beta u''&=\lambda u,\hspace{10pt} u(0)=u(1)=0
\end{align}
on $\Omega=(0,1)$. \autoref{int1d}(a) displays $\Omega$ with interface 
condition at $x=\frac{1}{3}$. $\Omega^+=(0,\frac{1}{3}), \Omega^{-}=(\frac{1}{3},1).$
$\Omega$ is divided into six spectral elements as displayed in \autoref{int1d}(b).\\
\subsubsection*{Problem 1:} Here we investigate \eqref{1d_int} with $\beta^{+}=1,\beta^{-}=4$  (\cite{zhang2022}). 
\autoref{res_1d} displays the obtained eigenvalues. First column represents the
polynomial order/degrees of freedom, second to seventh column represent first six
eigenvalues and last row represents the exact eigenvalues upto ten decimal places. $H^{1}$ error of first six eigenfunctions versus $W$ have been plotted in \autoref{er_1d_interface} in log scale.
\autoref{1d_int_results} presents slopes of the error plots.
\begin{table}[H]
\centering
\resizebox{0.98\linewidth}{!} {
    \begin{tabular}{D E E E E E E}
        \hline
        W/DOF & $\lambda_{1}$ & $\lambda_{2}$ & $\lambda_{3}$ & $\lambda_{4}$ & $\lambda_{5}$ & $\lambda_{6}$\\
        \hline
        4/30 &
        {\color{blue}22.2066}273051&
        {\color{blue}88.8}345949944&
        {\color{blue}199}.9555593711&
        {\color{blue}355}.4537100879&
        {\color{blue}559.1}055567547&
        827.7342917452\\
        6/42 &{\color{blue}22.20660990}42&
        {\color{blue}88.8264}401258&
        {\color{blue}199.859}5201113&
        {\color{blue}355.305}8220070&
        {\color{blue}555.1}766609487&
        {\color{blue}799}.6118359228\\
        8/54 &{\color{blue}22.20660990}18&
        {\color{blue}88.8264396}221&
        {\color{blue}199.8594891}408&
        {\color{blue}355.3057584}450&
        {\color{blue}555.1652}543130&
        {\color{blue}799.43}81770485				\\
        10/66 &{\color{blue}22.20660990}07&
        {\color{blue}88.8264396}181&
        {\color{blue}199.859489}0922&
        {\color{blue}355.30575843}35&
        {\color{blue}555.165247}6441&
        {\color{blue}799.437956}7570\\
        12/72 &{\color{blue}22.206609}8941&
        {\color{blue}88.82643960}43&
        {\color{blue}199.8594891}320&
        {\color{blue}355.30575843}37&
        {\color{blue}555.165247}4679&
        {\color{blue}799.437956}2572\\
        \hline
        {\color{blue}Exact}& {\color{blue}22.2066099024}&
        {\color{blue}88.8264396098}&
        {\color{blue}199.8594891221}&
        {\color{blue}355.3057584392}&
        {\color{blue}555.1652475613}&
        {\color{blue}799.4379564882}\\
        \hline
    \end{tabular}
    }
    \caption{Eigenvalues with $\beta^{+}=1, \beta^{-}=\frac{1}{4}$}
    \label{res_1d}
\end{table}
\begin{minipage}{\textwidth}
\centering
\begin{minipage}[!]{0.47\linewidth}
\centering
\begin{figure}[H]
\centering
\resizebox{0.85\linewidth}{!} {
\begin{tikzpicture}
				\begin{axis}[xlabel=$W$,
					ylabel=$\log_{10}\left(\|u_{i}-u_{i}^{W}\|_{H^{1}}\right)$,
					xmin=3, xmax=13,thick]
					\addplot[color=red,mark=square,thick] coordinates {
						(4,-4.2966)
						(6,-7.3076)
						(8,-9.8823)
						(10,-10.1950)
						(12,-9.8770)
					};
					\addplot[color=blue,mark=triangle,thick]coordinates{
						(4,-2.6382)
						(6,-5.0535)
						(8,-7.6884)
						(10,-9.1541)
						(12,-9.1541)
					};
					\addplot[color=green,mark=x,thick] coordinates {
						(4,-1.9050)
						(6, -3.9769)
						(8,-6.2605)
						(10,-8.7118)
						(12,-8.6238)
					};
					\addplot[color=teal,mark=o,thick] coordinates {
						(4,-1.6048)
						(6,-3.5088)
						(8,-5.6249)
						(10,-7.9135)
						(12,-8.5976)
					};
					\addplot[color=magenta,mark=diamond,thick] coordinates {
						(4,-0.7403)
						(6,-2.4360)
						(8,-4.2843)
						(10,-6.3070)
						(12,-8.1658)
					};
					\addplot[color=olive,mark=triangle*,thick] coordinates {
						(4,-0.3540)
						(6, -1.7442)
						(8,-3.4348)
						(10,-5.2986)
						(12,-7.2977)
					};
					\legend{$u_{1}$,$u_{2}$,$u_{3}$,$u_{4}$,$u_{5}$,$u_{6}$}
				\end{axis}
			\end{tikzpicture}
   }
\caption{$H^{1}$ error of eigenfunctions }
\label{er_1d_interface}
\end{figure}
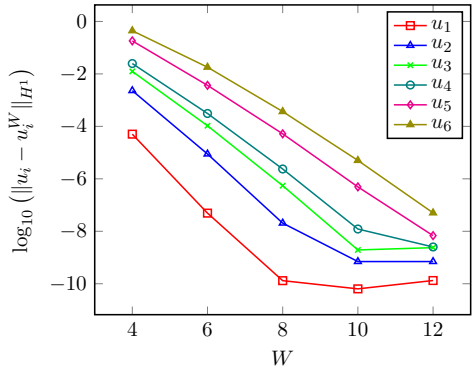
\end{minipage}
\begin{minipage}[!]{0.45\linewidth}
\begin{table}[H]
	\begin{center}
		\begin{tabular}{D C}
			\hline
   \hskip 0.1in
			plot & slope\\
			\hline 
			$u_1$ &  -0.7024\\
			$u_2$ &  -0.8643\\
			  $u_3$ &  -0.9086\\
			$u_4$ &  -0.9195\\
			$u_5$ &  -0.9361\\
	          $u_6$ &  -0.8721\\
   \hline 
		\end{tabular}
	\end{center}
	\caption{Slopes of error curves in \autoref{er_1d_interface}}
     \label{1d_int_results}
\end{table}
\end{minipage}
\end{minipage}
 \subsubsection*{Problem 2:}  Here we test \eqref{1d_int} with $\beta^{+}=1,\beta^{-}=1000$. \autoref{1:1000} presents first six eigenvalues using proposed method. We have neither exact nor reference eigenvalues with us, hence only obtained numerical values are reported upto ten decimal places.
    \begin{table}[H]
        \centering
        \resizebox{0.98\textwidth}{!}{
        \begin{tabular}{D E E E F G H}
             \hline
             W/DOF & $\lambda_{1}$ & $\lambda_{2}$ & $\lambda_{3}$ & $\lambda_{4}$ & $\lambda_{5}$ & $\lambda_{6}$\\
             \hline
             4/30 & 88.4793794235 & 354.0035048408 & 838.0772788461 &2127.1053071007 & 3721.4910570699 &5457.0488715691\\ 
            6/42 & 88.4674996888 & 353.8099254405& 796.0699555686& 1423.9806900782& 2252.7620572843& 3253.1749767088\\
            8/54 & 88.4675038184& 353.8098639823& 795.8139927784& 1414.0214004604& 2207.6175707870& 3172.5397961816\\
            10/66 & 88.4675110417& 353.8098191397& 795.8135025886& 1413.9758604099& 2207.0612570358& 3171.2822916881\\
            12/78 & 88.4675122041& 353.8098303336& 795.8142942182& 1413.9780204518& 2207.0577665966& 3171.2778934391\\
            \hline
        \end{tabular}
        }
        \caption{Eigenvalues with $\beta^{+}=1, \beta^{-}=1000$}
        \label{1:1000}
    \end{table}
\section{Conclusion and future work}
In this paper, we discuss a least-squares spectral element formulation for one-dimensional eigenvalue problems with interface conditions. Error estimates for eigenpairs have been discussed. Numerical experiments with different jumps have have been investigated. Eigenvalue problems with interface conditions on higher dimensional case will be presented in near future.

\end{document}